\documentclass{amsart}

\usepackage{amsfonts}
\usepackage{amsmath}
\usepackage{amssymb}
\usepackage{latexsym}
\usepackage{amscd}
\usepackage{amsthm}
\usepackage{subfigure}
\usepackage{graphicx}
\usepackage[all]{xy}
\usepackage{color}

\newtheorem{thm}{Theorem}[section]
\newtheorem{prop}[thm]{Proposition}
\newtheorem{cor}[thm]{Corollary}
\newtheorem{lem}[thm]{Lemma}

\theoremstyle{definition}

\newtheorem{rem}[thm]{Remark}

\newcommand{\Z}{\mathbb{Z}}
\newcommand{\M}{\mathcal{M}}
\newcommand{\T}{\mathcal{T}}
\newcommand{\I}{\mathcal{I}}
\newcommand{\G}{\Gamma_2}
\newcommand{\Y}{\mathcal{Y}}
\newcommand{\YLO}{\mathcal{Y}_{\mathrm{LO}}}

\author[N. Imoto]{Nao Imoto}

\author[R. Kobayashi]{Ryoma Kobayashi (Corresponding author)}
\address[R. Kobayashi]{
Department of General Education,\endgraf
National Institute of Technology, Ishikawa College,\endgraf
Tsubata, Ishikawa, 929-0392, Japan
}
\email{kobayashi\_ryoma@ishikawa-nct.ac.jp}



\thanks{The second author was supported by JSPS KAKENHI Grant Number JP19K14542 and JP22K13920.
Data sharing not applicable to this article as no datasets were generated or analysed during the current study.
The authors declare that they have no conflict of interest.}

\begin{document}

\title[On squares of Dehn twists about non-separating curves of $N_g$]
{On squares of Dehn twists about non-separating curves of a non-orientable closed surface}

\maketitle

\begin{abstract}
The level $2$ mapping class group of an orientable closed surface can be generated by squares of Dehn twists about non-separating curves (see \cite{Hu}).
On the other hand, the level $2$ mapping class group $\M_2(N)$ of a non-orientable closed surface $N$ can not be generated by only Dehn twists, and so it can not be generated by squares of Dehn twists about non-separating curves.
In this paper, we prove that the Dehn twist subgroup of $\M_2(N)$ can not be generated by squares of Dehn twists about non-separating curves either.
As an application, we give a finite generating set for the subgroup of $\M_2(N)$ generated by Dehn twists about separating curves and squares of Dehn twists about non-separating curves.
Moreover, we examine about actions on non-separating simple closed curves of $N$ by $\M_2(N)$.
\end{abstract}

\section{Introduction}

We first explain about the case of orientable surfaces.
For $g\geq0$, let $\Sigma_g$ denote an orientable closed surface of genus $g$, that is, $\Sigma_g$ is a connected sum of $g$ tori.
The \textit{mapping class group} $\M(\Sigma_g)$ of $\Sigma_g$ is the group which consists of isotopy classes of all orientation preserving diffeomorphisms of $\Sigma_g$.
The \textit{level $2$ mapping class group} $\M_2(\Sigma_g)$ of $\Sigma_g$ is the subgroup of $\M(\Sigma_g)$ which acts on $H_1(\Sigma_g;\Z/2\Z)$ trivially.
Let $\T^2(\Sigma_g)$ denote the subgroup of $\M(\Sigma_g)$ generated by squares of \textit{Dehn twists} about non-separating curves.
It is known that $\T^2(\Sigma_g)=\M_2(\Sigma_g)$ (see \cite{Hu}).
Let $\I(\Sigma_g)$ denote the subgroup of $\M(\Sigma_g)$ which acts on $H_1(\Sigma_g;\Z)$ trivially, called the \textit{Torelli group} of $\Sigma_g$.
It is clear that $\I(\Sigma_g)\subset\M_2(\Sigma_g)$, and hence we have that $\T^2(\Sigma_g)\I(\Sigma_g)=\M_2(\Sigma_g)$.

\begin{figure}[htbp]
\includegraphics{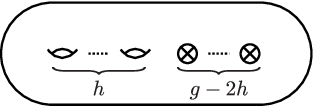}
\caption{A model of a non-orientable closed surface $N_g$.}\label{non-ori-surf}
\end{figure}

In this paper, we would like to consider the similar problem for the case of non-orientable surfaces, as follows.
For $g\geq1$, let $N_g$ denote a non-orientable closed surface of genus $g$, that is, $N_g$ is a connected sum of $\Sigma_h$ and $g-2h$ real projective planes for $\displaystyle0\leq{h}\leq\frac{g-1}{2}$.
In this paper, we describe $N_g$ as a surface obtained by attaching $g-2h$ M\"obius bands to a surface obtained by removing disjoint $g-2h$ disks from $\Sigma_h$, as shown in Figure~\ref{non-ori-surf}.
We call these attached M\"obius bands \textit{crosscaps}.
The \textit{mapping class group} $\M(N_g)$ of $N_g$ is the group which consists of isotopy classes of all diffeomorphisms of $N_g$.
The \textit{level $2$ mapping class group} $\M_2(N_g)$ of $N_g$ is the subgroup of $\M(N_g)$ which acts on $H_1(N_g;\Z/2\Z)$ trivially.
Let $\T^2(N_g)$ denote the subgroup of $\M(N_g)$ generated by squares of Dehn twists about non-separating curves.
Note that $\T^2(N_g)\subset\M_2(N_g)$.
$\M_2(N_g)$ can be normally generated by one \textit{crosscap slide} (see \cite{Sz2}), and a crosscap slide can not be described as any product of Dehn twists (see \cite{Li}).
Hence $\M_2(N_g)$ can not be generated by only Dehn twists, and so $\T^2(N_g)\neq\M_2(N_g)$.
We consider the subgroup $\T_2(N_g)$ of $\M_2(N_g)$ generated by all Dehn twists, called the \textit{twist subgroup} of $\M_2(N_g)$.
Note that a generating set for $\T_2(N_g)$ which does not consist of squares of Dehn twists about non-separating curves is known (see \cite{KO}).
In this paper, we prove that $\T_2(N_g)\neq\T^2(N_g)$ for $g\geq4$ either.

More precisely, we prove the following, where $\I(N_g)$ is the subgroup of $\M(N_g)$ which acts on $H_1(N_g;\Z)$ trivially, called the \textit{Torelli group} of $N_g$.

\begin{thm}\label{main-1}
For $g\geq3$, the quotient of $\M_2(N_g)$ by $\T^2(N_g)\I(N_g)$ is isomorphic to $\left(\Z/2\Z\right)^{\tbinom{g-1}{2}}$ if $g$ is odd and to $\left(\Z/2\Z\right)^{\tbinom{g-1}{2}+1}$ if $g$ is even.
\end{thm}

It is known that $\T_2(N_g)$ is an index $2$ subgroup of $\M_2(N_g)$ (see \cite{Li, KO}).
Since $\T_2(N_g)\supset\T^2(N_g)\I(N_g)$, we have the following.

\begin{cor}\label{main-1-cor}
For $g\geq3$, the quotient of $\T_2(N_g)$ by $\T^2(N_g)\I(N_g)$ is isomorphic to $\left(\Z/2\Z\right)^{\tbinom{g-1}{2}-1}$ if $g$ is odd and to $\left(\Z/2\Z\right)^{\tbinom{g-1}{2}}$ if $g$ is even.
\end{cor}

By this corollary, we can conclude that $\T_2(N_g)\neq\T^2(N_g)$ for $g\geq4$.
We remark that $\T_2(N_3)=\T^2(N_3)$ since $\I(N_3)$ is trivial.

Here is an outline of the paper.
In Section~\ref{pre}, we present the prior results necessary to prove Theorem~\ref{main-1}.
In Section~\ref{action}, we examine about actions on non-separating simple closed curves of $N_g$ by $\M_2(N_g)$.
In Section~\ref{main}, we prove Theorem~\ref{main-1}, using results of Sections~\ref{pre} and \ref{action}.
In Section~\ref{gen-T^2I}, we give a finite generating set for $\T^2(N_g)\I(N_g)$, as an application of Theorem~\ref{main-1}.
In Appendix~\ref{T^2I}, we consider the case where $g$ is odd of Theorem~\ref{main-1} and Corollary~\ref{main-1-cor}.

Through this paper, the product $gf$ of mapping classes $f$ and $g$ means that we apply $f$ first and then $g$, and we do not distinguish a simple closed curve from its isotopy class.
In addition, for simple closed curves $\alpha$ and $\beta$ of $N_g$, we denote $\alpha=\beta$ if $\alpha$ and $\beta$ are isotopic.

\section{Preliminaries}\label{pre}

\subsection{On Dehn twist and crosscap slide}\

\begin{figure}[htbp]
\subfigure[The Dehn twist $t_\alpha$ about $\alpha$.]{\includegraphics{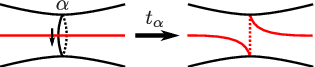}}~\subfigure[The crosscap slide $Y_{\mu,\alpha}$ about $\mu$ and $\alpha$.]{\includegraphics{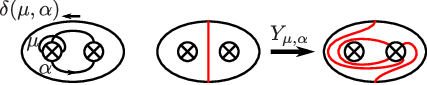}}
\caption{Descriptions of a Dehn twist and a crosscap slide.}\label{dehn-slide}
\end{figure}

First, we define a \textit{Dehn twist} and a \textit{crosscap slide} which are elements of $\M(N_g)$.
For a simple closed curve $\alpha$ of $N_g$, its regular neighborhood is either an annulus  or a M\"obius band.
We call $\alpha$ a \textit{two sided} or a \textit{one sided} simple closed curve respectively.
For a two sided simple closed curve $\alpha$, the \textit{Dehn twist} $t_\alpha$ about $\alpha$ is the isotopy class of the map acting as shown in Figure~\ref{dehn-slide}~(a).
The direction of the twist is indicated by an arrow written beside $\alpha$ as shown in Figure~\ref{dehn-slide}~(a).
For a one sided simple closed curve $\mu$ of $N_g$ and an oriented two sided simple closed curve $\alpha$ of $N_g$ such that $N_g\setminus\alpha$ is non-orientable when $g\geq3$ and that $\mu$ and $\alpha$ intersect transversely at only one point, the \textit{crosscap slide} $Y_{\mu,\alpha}$ about $\mu$ and $\alpha$ is the isotopy class of the map described by pushing the crosscap which is a regular neighborhood of $\mu$ once along $\alpha$, as shown in Figure~\ref{dehn-slide}~(b).

\begin{figure}[htbp]
\subfigure[]{\includegraphics{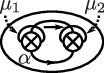}}\hspace{10mm}
\subfigure[]{\includegraphics{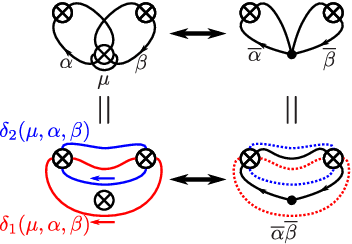}}
\caption{Relations on crosscap slides.}\label{crosscap-slide-rel}
\end{figure}

Next, we explain about some relations on Dehn twists and crosscap slides.
For details, for instance see \cite{Sz1}.
For a simple closed curve $\alpha$ which bounds a disk or a crosscap, we have that
$$t_\alpha=1.$$
For any Dehn twist $t_\alpha$, crosscap slide $Y_{\mu,\alpha}$ and $f\in\M(N_g)$, we have that
$$ft_\alpha{}f^{-1}=t_{f(\alpha)},~fY_{\mu,\alpha}f^{-1}=Y_{f(\mu),f(\alpha)},$$
where, the directions of the twist of $t_{f(\alpha)}$ and the pushing of $Y_{f(\mu),f(\alpha)}$ are induced from $f$ and the directions of the twist of $t_\alpha$ and the pushing of the $Y_{\mu,\alpha}$.
Let $\mu_1$, $\mu_2$ and $\alpha$ be simple closed curves as shown in Figure~\ref{crosscap-slide-rel}~(a).
Then we have that
$$Y_{\mu_2,\alpha}^{-1}Y_{\mu_1,\alpha}=Y_{\mu_2,\alpha}Y_{\mu_1,\alpha}^{-1}=t_\alpha^2.$$
For two crosscap slides $Y_{\mu,\alpha}$ and $Y_{\mu,\beta}$ such that $\overline{\alpha}\overline{\beta}$ is simple, where $\overline{\alpha}$ and $\overline{\beta}$ are oriented loops of $N_{g-1}$ obtained by collapsing a regular neighborhood of $\mu$ and $\overline{\alpha}\overline{\beta}$ is a composition loop of $\overline{\alpha}$ and $\overline{\beta}$ based at the collapsing point, let $\delta_1(\mu,\alpha,\beta)$ and $\delta_2(\mu,\alpha,\beta)$ be simple closed curves determined by $\mu$, $\alpha$ and $\beta$ as shown in Figure~\ref{crosscap-slide-rel}~(b).
Then we have that
$$Y_{\mu,\beta}Y_{\mu,\alpha}=t_{\delta_1(\mu,\alpha,\beta)}t_{\delta_2(\mu,\alpha,\beta)}.$$
In addition, for any crosscap slide $Y_{\mu,\alpha}$, we have that
$$Y_{\mu,\alpha}^2=t_{\delta(\mu,\alpha)},~Y_{\mu,\alpha}=Y_{\mu,\alpha^{-1}}^{-1},$$
where $\delta(\mu,\alpha)$ is a simple closed curve determined by $\mu$ and $\alpha$ as shown in Figure~\ref{dehn-slide}~(b).

\subsection{On $\M_2(N_g)$ and $\I(N_g)$}\label{M_2-I}\

\begin{figure}[htbp]
\includegraphics{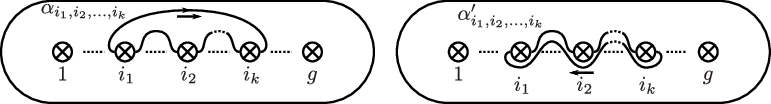}
\caption{Simple closed curves $\alpha_{i_1,i_2,\dots,i_k}$ and $\alpha^\prime_{i_1,i_2,\dots,i_k}$ for $1\leq{i_1}<i_2<\dots<i_k\leq{g}$.}\label{alpha}
\end{figure}

For $1\leq{i_1}<i_2<\dots<i_k\leq{g}$, let $\alpha_{i_1,i_2,\dots,i_k}$ be a simple closed curve of $N_g$ as shown in Figure~\ref{alpha}.
For $g\geq4$, Szepietowski \cite{Sz2} gave a finite generating set for $\M_2(N_g)$, and then Hirose-Sato \cite{HS} gave a minimal generating set for $\M_2(N_g)$ as follows.

\begin{thm}[\cite{HS}]\label{level2-gen}
For $g\geq4$, $\M_2(N_g)$ is minimally generated by
\begin{itemize}
\item	$Y_{\alpha_i,\alpha_{i,j}}$ for $1\leq{i<j}\leq{g}$,
\item	$Y_{\alpha_j,\alpha_{i,j}}$ for $1\leq{i<j}\leq{g-1}$, 
\item	$t_{\alpha_{1,j,k,l}}^2$ for $1<j<k<l\leq{g}$.
\end{itemize}
\end{thm}

For $g\leq3$, $\I(N_g)$ is trivial.
For $g\geq4$, a normal generating set for $\I(N_g)$ in $\M(N_g)$ is known as follows.

\begin{thm}[\cite{HK, Ko2}]\label{torelli-gen}
Let $\alpha$, $\beta_1$, $\beta_2$ and $\gamma$ be simple closed curves as shown in Figure~\ref{Torelli-gen}.
For $g\geq5$, $\I(N_g)$ is normally generated by $t_\alpha$ and $t_{\beta_1}t_{\beta_2}^{-1}$ in $\M(N_g)$.
$\I(N_4)$ is normally generated by $t_\alpha$, $t_{\beta_1}t_{\beta_2}^{-1}$ and $t_\gamma$ in $\M(N_4)$.

\begin{figure}[htbp]
\includegraphics{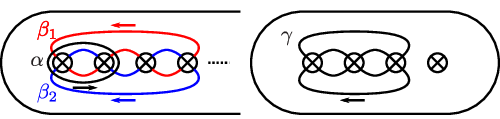}
\caption{Simple closed curves $\alpha$, $\beta_1$, $\beta_2$ and $\gamma$.}\label{Torelli-gen}
\end{figure}

\end{thm}

For $n\geq1$, let $\G(n)$  denote the kernel of the natural homomorphism $GL(n,\Z)\to{}GL(n,\Z/2\Z)$.
A finite presentation for $\G(n)$ was given by Fullarton \cite{Fu} and the second author \cite{Ko1} independently, and then Hirose and the second author \cite{HK} gave another presentation as follows.
Put $n=g-1$.
For $1\leq{i,j}\leq{g-1}$ with $i\neq{j}$, let $Y_{i,j}$ be the matrix whose $(i,j)$ entry is $2$, $(i,i)$ entry is $-1$, the other diagonal entries are $1$ and the other entries are $0$, and $Y_{i,g}$ the matrix whose $(i,i)$ entry is $-1$, the other diagonal entries are $1$ and the other entries are $0$.
$\G(g-1)$ has the following presentation.

\begin{prop}[\cite{HK}]\label{G_2-pre}
For $g\geq2$, $\G(g-1)$ admits a presentation with generators $Y_{i,j}$ for $1\leq{i}\leq{g-1}$ and $1\leq{j}\leq{g}$ with $i\neq{j}$.
The defining relators are
\begin{enumerate}
\item	$Y_{i,j}^2$,
\item	$(Y_{i,j}Y_{k,j})^2$, $(Y_{i,j}Y_{k,l})^2$,
\item	$(Y_{i,j}Y_{i,k}Y_{j,k})^2$, $(Y_{i,j}Y_{i,k}Y_{i,l})^2$,
\item	$(Y_{j,i}Y_{i,j}Y_{k,j}Y_{j,k}Y_{i,k}Y_{k,i})^2$,
\end{enumerate}
where indices $i$, $j$, $k$ and $l$ are distinct.
\end{prop}

We regard $\alpha_{i_1,i_2,\dots,i_k}$ as an element of $H_1(N_g;\Z)$.
Note that $\alpha_{i_1,i_2,\dots,i_k}=\alpha_{i_1}+\alpha_{i_2}+\cdots+\alpha_{i_k}$.
Then, as a $\Z$-module, we have presentations
\begin{eqnarray*}
H_1(N_g;\Z)/{\langle{\alpha_{1,2,\dots,g}}\rangle}&=&\langle{\alpha_1,\alpha_2,\dots,\alpha_g\mid\alpha_{1,2,\dots,g}=0}\rangle\\
&=&\langle{\alpha_1,\alpha_2,\dots,\alpha_{g-1}}\rangle.
\end{eqnarray*}
It is known that any automorphism of $H_1(N_g;\Z)$ preserves the homology class of $\alpha_{1,2,\dots,g}$ (see \cite{MP, HK}).
Hence $f\in\M(N_g)$ induces the automorphism $f_\ast$ of $H_1(N_g;\Z)/{\langle{\alpha_{1,2,\dots,g}}\rangle}$.
More precisely, for $f\in\M(N_g)$, when $f(\alpha_j)=\sum_{i=1}^gc_{ij}\alpha_i$, we define $f_\ast(\alpha_j)=\sum_{i=1}^{g-1}(c_{ij}-c_{gj})\alpha_i$.
We can regard an automorphism of $H_1(N_g;\Z)/{\langle{\alpha_{1,2,\dots,g}}\rangle}$ as a matrix in $GL(g-1,\Z)$.
Then we have that $(Y_{\alpha_i,\alpha_{i,j}})_\ast=Y_{i,j}$, $(Y_{\alpha_j,\alpha_{i,j}})_\ast=Y_{j,i}$ and $(Y_{\alpha_i,\alpha_{i,g}})_\ast=Y_{i,g}$ for $1\leq{i<j}\leq{g-1}$.
It is known that the quotient $\M_2(N_g)/{\I(N_g)}$ of $\M_2(N_g)$ by $\I(N_g)$ is isomorphic to $\G(g-1)$ for $g\geq3$ (see \cite{MP, HK}).

\section{Action on non-separating simple closed curves of $N_g$ by $\M_2(N_g)$}\label{action}

For a simple closed curve $\alpha$ of $N_g$, we denote by $[\alpha]$ its $\Z/2\Z$ coefficient first homology class.
In this section, we prove the following.

\begin{thm}\label{SCC/M_2}
For any non-separating simple closed curves $\alpha$ and $\beta$ of $N_g$, $[\alpha]=[\beta]\in{}H_1(N_g;\Z/2\Z)$ if and only if there is $\varphi\in\M_2(N_g)$ such that $\varphi(\alpha)=\beta$.
\end{thm}

\begin{proof}
We note that $H_1(N_g;\Z/2\Z)$ is isomorphic to $\left(\Z/2\Z\right)^{g}$.
Let $O_2(g)$ denote the orthogonal group on $\Z/2\Z$ of rank $g$.
The natural action on $H_1(N_g;\Z/2\Z)$ by $\M(N_g)$ induces the epimorphism $\M(N_g)\to{}O_2(g)$ (see \cite{GP}).
For $f\in\M(N_g)$, we denote by $f_\star\in{}O_2(g)$ the matrix corresponding to $f$.
It is known that $O_2(g)$ is generated by $(t_{\alpha_{i_1,i_2,\dots,i_k}})_\star$.
Since $\alpha_i$ and $\alpha_{i_1,i_2,\dots,i_k}$ intersect transversally at only one point if $i\in\{i_1,i_2,\dots,i_k\}$, and do not intersect if $i\notin\{i_1,i_2,\dots,i_k\}$, we have that
$$
(t_{\alpha_{i_1,i_2,\dots,i_k}})_\star([\alpha_i])=\left\{
\begin{array}{ll}
[\alpha_{i_1,i_2,\dots,i_k}]+[\alpha_i]&(i\in\{i_1,i_2,\dots,i_k\}),\\
{[\alpha_i]}&(i\notin\{i_1,i_2,\dots,i_k\}).
\end{array}
\right.
$$
We note that the kernel of $\M(N_g)\to{}O_2(g)$ is $\M_2(N_g)$.

Let $\alpha$ and $\beta$ be non-separating simple closed curves of $N_g$.
If there is $\varphi\in\M_2(N_g)$ such that $\varphi(\alpha)=\beta$, since $\varphi_\star=1$, we have that $[\alpha]=[\beta]$ clearly.
Hence we show the converse.
Suppose that $[\alpha]=[\beta]$.
There is $h\in\M(N_g)$ such that $h(\beta)$ is isotopic to either $\alpha_1$, $\alpha_{1,2}$ or $\alpha_{1,2,\dots,g}$.
We construct $f\in\M_2(N_g)$ satisfying $f(h(\alpha))=h(\beta)$.

First, suppose that $h(\beta)=\alpha_1$.
Note that $[h(\alpha)]=[h(\beta)]=[\alpha_1]$.
For $f_1\in\M(N_g)$ satisfying $f_1(h(\alpha))=h(\beta)$, we see that $(f_1)_\star([\alpha_1])=(f_1)_\star([h(\alpha)])=[f_1(h(\alpha))]=[h(\beta)]=[\alpha_1]$.
On the other hand, for $2\leq{i}\leq{g}$, since the mod $2$ intersection number of $(f_1)_\star([\alpha_1])$ and $(f_1)_\star([\alpha_i])$ is $0$, we can describe $\displaystyle(f_1)_\star([\alpha_i])=\sum_{j=2}^gc_{ji}[\alpha_j]$, where $c_{ji}=0$ or $1$.
Note that $(f_1)_\star$ can be described as a matrix
$$
\begin{pmatrix}
1&0&\cdots&0\\
0&c_{22}&\cdots&c_{2g}\\
\vdots&\vdots&&\vdots\\
0&c_{g2}&\cdots&c_{gg}
\end{pmatrix}.
$$
Then there is $f_2\in\M(N_g)$ such that
\begin{itemize}
\item	$(f_2)_\star=(f_1)_\star^{-1}$,
\item	$f_2$ is a product of $t_{\alpha_{i_1,i_2,\dots,i_k}}$'s for some even indices $i_1$, $i_2$, $\dots$, $i_k$ with $i_1\geq2$.
\end{itemize}
Let $f=f_2f_1$.
We note that $f_\star=(f_2)_\star(f_1)_\star=1$, that is, $f$ is in $\M_2(N_g)$.
Since $t_{\alpha_{i_1,i_2,\dots,i_k}}(\alpha_1)=\alpha_1$ if $i_1\geq2$, we have that
$$f(h(\alpha))=f_2(h(\beta))=f_2(\alpha_1)=\alpha_1=h(\beta).$$

Next, suppose that $h(\beta)=\alpha_{1,2}$.
Note that $[h(\alpha)]=[h(\beta)]=[\alpha_1]+[\alpha_2]$.
For $f_1\in\M(N_g)$ satisfying $f_1(h(\alpha))=h(\beta)$, we see that $(f_1)_\star([\alpha_1]+[\alpha_2])=(f_1)_\star([h(\alpha)])=[f_1(h(\alpha))]=[h(\beta)]=[\alpha_1]+[\alpha_2]$.
Put $\displaystyle(f_1)_\star([\alpha_i])=\sum_{j=1}^gc_{ji}[\alpha_j]$ for $1\leq{i}\leq{g}$, where $c_{ji}=0$ or $1$.
By the definition of $O_2(g)$, we have that there are an odd number of $c_{j1}$'s (resp. $c_{j2}$'s) for which $c_{j1}=1$ (resp. $c_{j2}=1$), and that there are an even number of $j$'s for which $c_{j1}=c_{j2}=1$.
In addition, we see that $\displaystyle(f_1)_\star([\alpha_1]+[\alpha_2])=\sum_{j=1}^g(c_{j1}+c_{j2})[\alpha_j]=[\alpha_1]+[\alpha_2]$, and hence $c_{11}\neq{}c_{12}$, $c_{21}\neq{}c_{22}$ and $c_{j1}=c_{j2}$ for $j\geq3$.
Therefore, we have that
$
\begin{pmatrix}
c_{11}&c_{12}\\
c_{21}&c_{22}
\end{pmatrix}
=
\begin{pmatrix}
1&0\\
0&1
\end{pmatrix}
$
or
$
\begin{pmatrix}
0&1\\
1&0
\end{pmatrix}
$.
We can rewrite $\displaystyle(f_1)_\star([\alpha_1])=[\alpha_i]+\sum_{s=1}^l[\alpha_{j_s}]$ for some even indices $3\leq{j_1<j_2<\dots<j_l}\leq{g}$, where $i=1$ or $2$.
Then we calculate
\begin{eqnarray*}
(t_{\alpha_{1,2,j_1,\dots,j_l}})_\star(f_1)_\star([\alpha_1])
&=&(t_{\alpha_{1,2,j_1,\dots,j_l}})_\star([\alpha_i])+\sum_{s=1}^l(t_{\alpha_{1,2,j_1,\dots,j_l}})_\star([\alpha_{j_s}])\\
&=&[\alpha_{1,2,j_1,\dots,j_l}]+[\alpha_i]+\sum_{s=1}^l([\alpha_{1,2,j_1,\dots,j_l}]+[\alpha_{j_s}])\\
&=&(l+1)[\alpha_{1,2,j_1,\dots,j_l}]+\left([\alpha_i]+\sum_{s=1}^l[\alpha_{j_s}]\right)\\
&=&(l+2)[\alpha_{1,2,j_1,\dots,j_l}]-[\alpha_j]\\
&=&[\alpha_j],
\end{eqnarray*}
where $(i,j)=(1,2)$ or $(2,1)$.
Similarly $(t_{\alpha_{1,2,j_1,\dots,j_l}})_\star(f_1)_\star([\alpha_2])=[\alpha_i]$.
Note that $(t_{\alpha_{1,2,j_1,\dots,j_l}})_\star(f_1)_\star$ can be described as a matrix
$$
\begin{pmatrix}
1&0&0&\cdots&0\\
0&1&0&\cdots&0\\
0&0&c_{33}^\prime&\cdots&c_{3g}^\prime\\
\vdots&\vdots&\vdots&&\vdots\\
0&0&c_{g3}^\prime&\cdots&c_{gg}^\prime
\end{pmatrix}
~\textrm{or}~
\begin{pmatrix}
0&1&0&\cdots&0\\
1&0&0&\cdots&0\\
0&0&c_{33}^\prime&\cdots&c_{3g}^\prime\\
\vdots&\vdots&\vdots&&\vdots\\
0&0&c_{g3}^\prime&\cdots&c_{gg}^\prime
\end{pmatrix}.
$$
Then there is $f_2\in\M(N_g)$ such that
\begin{itemize}
\item	$(f_2)_\star=\{(t_{\alpha_{1,2,j_1,\dots,j_l}})_\star(f_1)_\star\}^{-1}$,
\item	$f_2$ is a product of $t_{\alpha_{i_1,i_2,\dots,i_k}}$'s (and $t_{\alpha_{1,2}}$) for some even indices $i_1$, $i_2$, $\dots$, $i_k$ with $i_1\geq3$.
\end{itemize}
Let $f=f_2t_{\alpha_{1,2,j_1,\dots,j_l}}f_1$.
We note that $f_\star=(f_2)_\star\{(t_{\alpha_{1,2,j_1,\dots,j_l}})_\star(f_1)_\star\}=1$, that is, $f$ is in $\M_2(N_g)$.
Since $t_{\alpha_{1,2,j_1,\dots,j_l}}(\alpha_{1,2})=\alpha_{1,2}$, $t_{\alpha_{1,2}}(\alpha_{1,2})=\alpha_{1,2}$ and $t_{\alpha_{i_1,i_2,\dots,i_k}}(\alpha_{1,2})=\alpha_{1,2}$ if $i_1\geq3$, we have that
$$f(h(\alpha))=f_2t_{\alpha_{1,2,j_1,\dots,j_l}}(h(\beta))=f_2t_{\alpha_{1,2,j_1,\dots,j_l}}(\alpha_{1,2})=f_2(\alpha_{1,2})=\alpha_{1,2}=h(\beta).$$

Finally, suppose that $h(\beta)=\alpha_{1,2,\dots,g}$.
Note that $[h(\alpha)]=[h(\beta)]=[\alpha_1]+[\alpha_2]+\cdots+[\alpha_g]$.
For $f_1\in\M(N_g)$ satisfying $f_1(h(\alpha))=h(\beta)$, there is $f_2\in\M(N_g)$ such that
\begin{itemize}
\item	$(f_2)_\star=(f_1)_\star^{-1}$,
\item	$f_2$ is a product of $t_{\alpha_{i_1,i_2,\dots,i_k}}$'s for some even indices $i_1$, $i_2$, $\dots$, $i_k$.
\end{itemize}
Let $f=f_2f_1$.
We note that $f_\star=(f_2)_\star(f_1)_\star=1$, that is, $f$ is in $\M_2(N_g)$.
Since $t_{\alpha_{i_1,i_2,\dots,i_k}}(\alpha_{1,2,\dots,g})=\alpha_{1,2,\dots,g}$, we have that
$$f(h(\alpha))=f_2(h(\beta))=f_2(\alpha_{1,2,\dots,g})=\alpha_{1,2,\dots,g}=h(\beta).$$

Therefore $\varphi=h^{-1}fh$ is in $\M_2(N_g)$ and satisfies $\varphi(\alpha)=\beta$ if $[\alpha]=[\beta]$.
Thus the proof is completed.
\end{proof}

\section{Proof of Theorem~\ref{main-1}}\label{main}

In this section, we prove Theorem~\ref{main-1}.
First we prepare some things.

For $1\leq{i}\leq{g-1}$, let $Y_{g,i}=(Y_{\alpha_g,\alpha_{i,g}})_\ast\in\G(g-1)$, and we can check that
$$Y_{g,i}=(Y_{1,i}Y_{1,g})(Y_{2,i}Y_{2,g})\cdots(Y_{i-1,i}Y_{i-1,g})\cdot(Y_{i+1,i}Y_{i+1,g})\cdots(Y_{g-1,i}Y_{g-1,g})Y_{i,g}.$$
Hence, from Proposition~\ref{G_2-pre} we have the following.

\begin{cor}\label{G_2-pre2}
For $g\geq3$, $\G(g-1)$ admits a presentation with generators $Y_{i,j}$ for $1\leq{i,j}\leq{g}$ with $i\neq{j}$.
The defining relators are $(1)$-$(4)$ of Proposition~\ref{G_2-pre} and
\begin{enumerate}
\item[(5)]	$(Y_{1,i}Y_{1,g})(Y_{2,i}Y_{2,g})\cdots(Y_{i-1,i}Y_{i-1,g})\cdot(Y_{i+1,i}Y_{i+1,g})\cdots(Y_{g-1,i}Y_{g-1,g})Y_{i,g}Y_{g,i}^{-1}$ for $1\leq{i}\leq{g-1}$.
\end{enumerate}
\end{cor}

\begin{lem}\label{G_2-rel}
For $g\geq3$, in $\G(g-1)$, the relation $Y_{i,j}Y_{g,j}=Y_{g,j}Y_{i,j}$ holds for $1\leq{i}\leq{g-1}$.
\end{lem}

\begin{proof}
By relators of $\G(g-1)$, we have relations
\begin{eqnarray*}
Y_{i,j}(Y_{k,j}Y_{k,g})&=&(Y_{k,j}Y_{k,g})Y_{i,j},\\
Y_{i,j}(Y_{i,j}Y_{i,g})&=&(Y_{i,j}Y_{i,g})(Y_{i,g}Y_{i,j}Y_{i,g}),\\
(Y_{i,g}Y_{i,j}Y_{i,g})(Y_{k,j}Y_{k,g})&=&(Y_{k,j}Y_{k,g})(Y_{i,g}Y_{i,j}Y_{i,g})
\end{eqnarray*}
for $k\neq{i}$.
Hence we calculate
\begin{eqnarray*}
Y_{i,j}Y_{g,j}
&=&Y_{i,j}(Y_{1,j}Y_{1,g})\cdots(Y_{j-1,j}Y_{j-1,g})\cdot(Y_{j+1,j}Y_{j+1,g})\cdots(Y_{g-1,j}Y_{g-1,g})Y_{j,g}\\
&=&(Y_{1,j}Y_{1,g})\cdots(Y_{j-1,j}Y_{j-1,g})\cdot(Y_{j+1,j}Y_{j+1,g})\cdots(Y_{g-1,j}Y_{g-1,g})(Y_{i,g}Y_{i,j}\underset{(2)}{\underline{Y_{i,g})Y_{j,g}}}\\
&=&(Y_{1,j}Y_{1,g})\cdots(Y_{j-1,j}Y_{j-1,g})\cdot(Y_{j+1,j}Y_{j+1,g})\cdots(Y_{g-1,j}Y_{g-1,g})\underset{(3)}{\underline{Y_{i,g}Y_{i,j}Y_{j,g}}}Y_{i,g}\\
&=&(Y_{1,j}Y_{1,g})\cdots(Y_{j-1,j}Y_{j-1,g})\cdot(Y_{j+1,j}Y_{j+1,g})\cdots(Y_{g-1,j}Y_{g-1,g})Y_{j,g}Y_{i,j}\underset{(1)}{\underline{Y_{i,g}Y_{i,g}}}\\
&=&(Y_{1,j}Y_{1,g})\cdots(Y_{j-1,j}Y_{j-1,g})\cdot(Y_{j+1,j}Y_{j+1,g})\cdots(Y_{g-1,j}Y_{g-1,g})Y_{j,g}Y_{i,j}\\
&=&Y_{g,j}Y_{i,j}.
\end{eqnarray*}
Thus we obtain the claim.
\end{proof}

Let $N$ be the normal subgroup of $\G(g-1)$ normally generated by $Y_{j,i}^{-1}Y_{i,j}$ for $1\leq{i,j}\leq{g}$.

\begin{lem}\label{G/N-1}
For $g\geq3$, the quotient $\G(g-1)/{N}$ of $\G(g-1)$ by $N$ admits a presentation with generators $Y_{i,j}$ for $1\leq{i<j}\leq{g}$.
The defining relators are
\begin{enumerate}
\item[$\overline{(1)}$]	$Y_{i,j}^2$,
\item[$\overline{(2)}$]	$(Y_{i,j}Y_{k,j})^2$, $(Y_{i,j}Y_{i,l})^2$, $(Y_{i,j}Y_{j,l})^2$, $(Y_{i,j}Y_{k,l})^2$,
\item[$\overline{(5)}$]	$(Y_{1,i}Y_{1,g})\cdots(Y_{i-1,i}Y_{i-1,g})\cdot(Y_{i,i+1}Y_{i+1,g})\cdots(Y_{i,g-1}Y_{g-1,g})$ for $2\leq{i}\leq{g-1}$, and $Y_{1,g}Y_{2,g}\cdots{}Y_{g-1,g}$ if $g$ is odd.
\end{enumerate}
\end{lem}

\begin{proof}
Since we have that $Y_{j,i}=Y_{i,j}$ for $1\leq{i<j}\leq{g}$ in $\G(g-1)/{N}$, by the Tietze transformation, $\G(g-1)/{N}$ is generated by $Y_{i,j}$ for $1\leq{i<j}\leq{g}$.

By the relator~(1) of $\G(g-1)$, we have the relator~$\overline{(1)}$.
By the relator $(Y_{i,j}Y_{k,j})^2$ of $\G(g-1)$, we have relators
\begin{itemize}
\item	$(Y_{i,j}Y_{k,j})^2$ if $1\leq{i,k}<j\leq{g}$,
\item	$(Y_{j,i}Y_{j,k})^2$ if $1\leq{j}<{i,k}\leq{g-1}$,
\item	$(Y_{i,j}Y_{j,k})^2$ if $1\leq{i<j<k}\leq{g-1}$.
\end{itemize}
In addition, by Lemma~\ref{G_2-rel}, we have relators
\begin{itemize}
\item	$(Y_{j,i}Y_{j,g})^2$ if $i>j$,
\item	$(Y_{i,j}Y_{j,g})^2$ if $i<j$.
\end{itemize}
Moreover, by the relator $(Y_{i,j}Y_{k,l})^2$ of $\G(g-1)$, we have the relator $(Y_{i,j}Y_{k,l})^2$.
Hence we have the relator~$\overline{(2)}$.
By the relators~$\overline{(1)}$ and $\overline{(2)}$, the relators~(3) and (4) are trivial relators.
By the relator~(5) of $\G(g-1)$, we have the relator
$$(Y_{1,i}Y_{1,g})\cdots(Y_{i-1,i}Y_{i-1,g})\cdot(Y_{i,i+1}Y_{i+1,g})\cdots(Y_{i,g-1}Y_{g-1,g}),$$
and hence by the relators~$\overline{(1)}$ and $\overline{(2)}$,
\begin{eqnarray*}
Y_{1,i}&=&Y_{1,g}(Y_{2,i}Y_{2,g})\cdots(Y_{i-1,i}Y_{i-1,g})\cdot(Y_{i,i+1}Y_{i+1,g})\cdots(Y_{i,g-1}Y_{g-1,g})\\
&=&(Y_{1,g}\cdots{}Y_{i-1,g}Y_{i+1,g}\cdots{}Y_{g-1,g})(Y_{2,i}\cdots{}Y_{i-1,i}Y_{i,i+1}\cdots{}Y_{i,g-1}).
\end{eqnarray*}
for $2\leq{i}\leq{g-1}$.
By the relators~$\overline{(1)}$ and $\overline{(2)}$, the relator~(5) of $\G(g-1)$ for $i=1$ can be transformed as
\begin{eqnarray*}
(Y_{1,2}Y_{2,g})\cdots(Y_{1,g-1}Y_{g-1,g})
&=&(Y_{1,2}\cdots{}Y_{1,g-1})(Y_{2,g}\cdots{}Y_{g-1,g})\\
&=&(Y_{1,g}\cdot{}Y_{3,g}\cdots{}Y_{g-1,g})(Y_{2,3}\cdots{}Y_{2,g-1})\\
&&(Y_{1,g}Y_{2,g}\cdot{}Y_{4,g}\cdots{}Y_{g-1,g})(Y_{2,3}\cdot{}Y_{3,4}\cdots{}Y_{3,g-1})\\
&&(Y_{1,g}Y_{2,g}Y_{3,g}\cdot{}Y_{5,g}\cdots{}Y_{g-1,g})(Y_{2,4}Y_{3,4}\cdot{}Y_{4,5}\cdots{}Y_{4,g-1})\\
&&\vdots\\
&&(Y_{1,g}\cdots{}Y_{g-3,g}\cdot{}Y_{g-1,g})(Y_{2,g-2}\cdots{}Y_{g-3,g-2}\cdot{}Y_{g-2,g-1})\\
&&(Y_{1,g}\cdots{}Y_{g-2,g})(Y_{2,g-1}\cdots{}Y_{g-2,g-1})\\
&&(Y_{2,g}\cdots{}Y_{g-1,g})\\
&=&(Y_{1,g}\cdots{}Y_{g-1,g})^{g-2}\prod_{2\leq{i<j}\leq{g-1}}Y_{i,j}^2\\
&=&
\left\{
\begin{array}{ll}
Y_{1,g}\cdots{}Y_{g-1,g}&\textrm{if}~g~\textrm{is odd},\\
1&\textrm{if}~g~\textrm{is even}.\\
\end{array}
\right.
\end{eqnarray*}
Hence we have the relator~$\overline{(5)}$.

Thus we obtain the claim.
\end{proof}

By the relators~$\overline{(1)}$ and $\overline{(2)}$, it follows that $\G(g-1)/{N}$ is a $\Z/2\Z$-module.
By the relator~$\overline{(5)}$, generators $Y_{1,i}$ are not needed, for $2\leq{i}\leq{g}$ if $g$ is odd and for $2\leq{i}\leq{g-1}$ if $g$ is even.
Therefore we have the following.

\begin{cor}\label{G/N-2}
For $g\geq3$, the quotient $\G(g-1)/{N}$ of $\G(g-1)$ by $N$ is isomorphic to $\left(\Z/2\Z\right)^{\tbinom{g-1}{2}}$ if $g$ is odd and to $\left(\Z/2\Z\right)^{\tbinom{g-1}{2}+1}$ if $g$ is even.
\end{cor}

\begin{lem}\label{norm-gen-T^2}
In $\M_2(N_g)$, $\T^2(N_g)$ is normally generated by $t_{\alpha_{i_1,i_2,\dots,i_k}}^2$ for even indices $1\leq{i_1<i_2<\cdots<i_k\leq{g}}$.
\end{lem}

\begin{proof}
For any non-separating simple closed curve $\alpha$ of $N_g$, there are indices $1\leq{i_1<i_2<\cdots<i_k\leq{g}}$ such that $[\alpha]=[\alpha_{i_1,i_2,\dots,i_k}]\in{}H_2(N_g;\Z/2\Z)$, and then there is $f\in\M_2(N_g)$ such that $f(\alpha_{i_1,i_2,\dots,i_k})=\alpha$, by Theorem~\ref{SCC/M_2}.
Hence for any $t_\alpha^2\in\T^2(N_g)$, there is $f\in\M_2(N_g)$ such that $t_\alpha^2=t_{f(\alpha_{i_1,i_2,\dots,i_k})}^2=ft_{\alpha_{i_1,i_2,\dots,i_k}}^2f^{-1}$.
Thus we obtain the claim.
\end{proof}

\begin{lem}\label{chain-T^2}
Let $d_1$, $d_2$ and $d_3$ be simple closed curves on a surface $S$ as shown in Figure~\ref{chain}.
Then $t_{d_1}t_{d_2}$ and $t_{d_3}$ are in $\T^2(S)$.

\begin{figure}[htbp]
\includegraphics{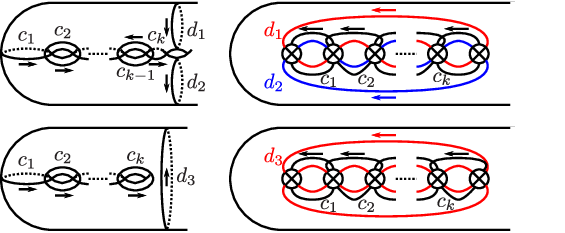}
\caption{Simple closed curves $c_1$, $c_2$, $\dots$, $c_k$, $d_1$, $d_2$ and $d_3$.}\label{chain}
\end{figure}

\end{lem}

\begin{proof}
Let $c_i$ be a simple closed curve of $S$ as shown in Figure~\ref{chain}, for $1\leq{i}\leq{k}$.
It is known that
\begin{itemize}
\item	$t_{d_1}t_{d_2}=(t_{c_1}t_{c_2}\cdots{}t_{c_k})^{k+1}$,
\item	$t_{d_3}=(t_{c_1}t_{c_2}\cdots{}t_{c_k})^{2(k+1)}$,
\end{itemize}
called the $k$-chain relation.
When $k=1$, it is clear that $t_{d_1}t_{d_2}$ is in $\T^2(S)$.
When $k\geq2$, by the commutativity relation $t_{c_i}t_{c_j}=t_{c_j}t_{c_i}$ with $|i-j|\geq2$ and the braid relation $t_{c_i}t_{c_{i+1}}t_{c_i}=t_{c_{i+1}}t_{c_i}t_{c_{i+1}}$, we calculate
\begin{eqnarray*}
(t_{c_1}t_{c_2}\cdots{}t_{c_k})^{k+1}
&=&(t_{c_1}t_{c_2}\cdots{}t_{c_{k-1}})^kt_{c_k}t_{c_{k-1}}\cdots{}t_{c_2}t_{c_1}^2t_{c_2}\cdots{}t_{c_{k-1}}t_{c_k}\\
&=&(t_{c_1}t_{c_2}\cdots{}t_{c_{k-1}})^k\\
&&t_{c_k}^2(t_{c_k}^{-1}t_{c_{k-1}}^2t_{c_k})\cdots(t_{c_k}^{-1}t_{c_{k-1}}^{-1}\cdots{}t_{c_2}^{-1}t_{c_1}^2t_{c_2}\cdots{}t_{c_{k-1}}t_{c_k}).
\end{eqnarray*}
By induction on $k$, we have that $(t_{c_1}t_{c_2}\cdots{}t_{c_k})^{k+1}$ is in $\T^2(S)$, and so are $t_{d_1}t_{d_2}$ and $t_{d_3}$.
Thus we obtain the claim.
\end{proof}

We now prove Theorem~\ref{main-1}.

\begin{proof}[Proof of Theorem~\ref{main-1}]
By Corollary~\ref{G/N-2}, it suffices to show that the quotient of $\M_2(N_g)$ by $\T^2(N_g)\I(N_g)$ is isomorphic to $\G(g-1)/{N}$.
Let $\phi:\M_2(N_g)\to\G(g-1)$ be the homomorphism defined in Subsection~\ref{M_2-I}, that is, $\phi(f)=f_\ast$.
Since $\ker\phi$ is $\I(N_g)$, we show that $\phi(\T^2(N_g))$ is $N$.

For $1\leq{i<j}\leq{g}$, we have that $(Y_{\alpha_j,\alpha_{i,j}}^{-1})_\ast(Y_{\alpha_i,\alpha_{i,j}})_\ast=Y_{j,i}^{-1}Y_{i,j}\in{N}$.
In addition, since the equality $t_{\alpha_{i,j}}^2=Y_{\alpha_j,\alpha_{i,j}}^{-1}Y_{\alpha_i,\alpha_{i,j}}$ holds, we have that $\phi(t_{\alpha_{i,j}}^2)=Y_{j,i}^{-1}Y_{i,j}$.
Therefore $\phi(\T^2(N_g))\supset{N}$.

In order to prove $\phi(\T^2(N_g))\subset{N}$, it suffices to show that $\phi(t_{\alpha_{i_1,i_2,\dots,i_k}}^2)$ is in $N$, by Lemma~\ref{norm-gen-T^2}.
Let $\alpha^\prime_{i_1,i_2,\dots,i_k}$ be a simple closed curve as shown in Figure~\ref{alpha}.
Since $t_{\alpha_{i_1,i_2,\dots,i_k}}^{-1}t_{\alpha^\prime_{i_1,i_2,\dots,i_k}}$ is in $\I(N_g)$ and hence its image by $\phi$ is $1$, we have that
$$\phi(t_{\alpha_{i_1,i_2,\dots,i_k}}^2)=\phi(t_{\alpha_{i_1,i_2,\dots,i_k}}^2\cdot{}t_{\alpha_{i_1,i_2,\dots,i_k}}^{-1}t_{\alpha^\prime_{i_1,i_2,\dots,i_k}})=\phi(t_{\alpha_{i_1,i_2,\dots,i_k}}t_{\alpha^\prime_{i_1,i_2,\dots,i_k}}).$$
By the proof of Lemma~\ref{chain-T^2}, $t_{\alpha_{i_1,i_2,\dots,i_k}}t_{\alpha^\prime_{i_1,i_2,\dots,i_k}}$ is a product of $T_{s,t}$ for $1\leq{s<t}\leq{k}$, where
$$T_{s,t}=(t_{\alpha_{i_{s+1},i_{s+2}}}t_{\alpha_{i_{s+2},i_{s+3}}}\cdots{}t_{\alpha_{i_{t-1},i_t}})^{-1}t_{\alpha_{i_s,i_{s+1}}}^2(t_{\alpha_{i_{s+1},i_{s+2}}}t_{\alpha_{i_{s+2},i_{s+3}}}\cdots{}t_{\alpha_{i_{t-1},i_t}}).$$
Note that $T_{s,t}$ is the square of the Dehn twist about a simple closed curve as shown in Figure~\ref{st}.
Similar to the proof of Lemma~\ref{norm-gen-T^2}, there is $f\in\M_2(N_g)$ such that $T_{s,t}=ft_{\alpha_{i_s,i_t}}^2f^{-1}$.
Since $\phi(t_{\alpha_{i_s,i_t}}^2)=Y_{i_t,i_s}^{-1}Y_{i_s,i_t}$, it follows that $\phi(t_{\alpha_{i_1,i_2,\dots,i_k}}t_{\alpha^\prime_{i_1,i_2,\dots,i_k}})$ is in $N$.
Therefore	 $\phi(\T^2(N_g))\subset{N}$.

\begin{figure}[htbp]
\includegraphics{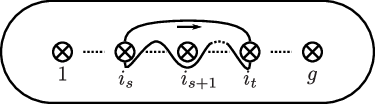}
\caption{The loop $(t_{\alpha_{i_{s+1},i_{s+2}}}t_{\alpha_{i_{s+2},i_{s+3}}}\cdots{}t_{\alpha_{i_{t-1},i_t}})^{-1}(\alpha_{i_s,i_{s+1}})$.}\label{st}
\end{figure}

Thus we complete the proof.
\end{proof}

\section{A finite generating set for $\T^2(N_g)\I(N_g)$}\label{gen-T^2I}

By Theorem~\ref{main-1}, $\T^2(N_g)\I(N_g)$ is a finite index subgroup of $\M_2(N_g)$.
In addition, by Theorem~\ref{level2-gen}, $\M_2(N_g)$ is finitely generated.
Hence $\T^2(N_g)\I(N_g)$ is also finitely generated.
In this section, we give a finite generating set for $\T^2(N_g)\I(N_g)$ by the Reidemeister-Schreier method.

\begin{figure}[htbp]
\includegraphics{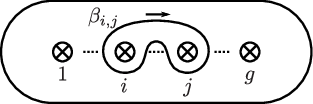}
\caption{A simple closed curve $\beta_{i,j}$ for $1\leq{i<j}\leq{g}$.}\label{beta_ij}
\end{figure}

For $1\leq{i<j}\leq{g}$, let $\beta_{i,j}$ be a simple closed curve as shown in Figure~\ref{beta_ij}.
For $(i,j)$ with $1\leq{i<j}\leq{g}$, we consider the lexicographic order $<$, that is, $(i,j)<(k,l)$ if and only if $i<k$, or $i=k$ and $j<l$.
Let $\Y$ be the set consisting of
$$Y_{\alpha_{i_1},\alpha_{i_1,j_1}}Y_{\alpha_{i_2},\alpha_{i_2,j_2}}\cdots{}Y_{\alpha_{i_m},\alpha_{i_m,j_m}}$$
for $2\leq{i_s<j_s}\leq{g}$ if $g$ is odd and $(i_s,j_s)=(1,g)$ or $2\leq{i_s<j_s}\leq{g}$ if $g$ is even, where $s\neq{}t\iff(i_s,j_s)\neq(i_t,j_t)$, and $\YLO$ the subset of $\Y$ consisting of
$$Y_{\alpha_{i_1},\alpha_{i_1,j_1}}Y_{\alpha_{i_2},\alpha_{i_2,j_2}}\cdots{}Y_{\alpha_{i_m},\alpha_{i_m,j_m}}\in\Y,$$
where $s<t\iff(i_s,j_s)<(i_t,j_t)$.

The main result of this section is as follows.

\begin{thm}\label{gen-T^2I-2}
$\T^2(N_g)\I(N_g)$ is generated by
\begin{enumerate}
\item	$ft_{\alpha_{i,j}}^2f^{-1}$ for $1\leq{i<j}\leq{g}$ and $f\in\YLO$,
\item	$ft_{\beta_{i,j}}f^{-1}$ for $1\leq{i<j}\leq{g}$ and $f\in\YLO$,
\item	$ft_{\alpha_{1,j,k,l}}^2f^{-1}$ for $1<j<k<l\leq{g}$ and $f\in\YLO$.
\end{enumerate}
\end{thm}

In order to prove this theorem, we first show the following lemma.

\begin{lem}\label{gen-T^2I-1}
$\T^2(N_g)\I(N_g)$ is generated by $(1)$-$(3)$ of Theorem~\ref{gen-T^2I-2} and
\begin{enumerate}
\item[(4)]	$f[Y_{\alpha_{i},\alpha_{i,j}},Y_{\alpha_{k},\alpha_{k,l}}]f^{-1}$ for $1\leq{i<j}\leq{g}$, $1\leq{k<l}\leq{g}$ and $f\in\YLO$.
\end{enumerate}
\end{lem}

In addition, in order to prove this lemma, we use the following lemma.

\begin{lem}\label{[a,b]}
The following holds.
\begin{eqnarray*}
[a,b_1b_2\cdots{}b_n]
&=&
[a,b_1]\left(b_1[a,b_2]b_1^{-1}\right)\left(b_1b_2[a,b_3]b_2^{-1}b_1^{-1}\right)\\
&&
\cdots\left(b_1b_2\cdots{}b_{n-1}[a,b_n]b_{n-1}^{-1}\cdots{}b_2^{-1}b_1^{-1}\right)
\end{eqnarray*}
\end{lem}

\begin{proof}
We calculate
\begin{eqnarray*}
&&[a,b_1]\left(b_1[a,b_2]b_1^{-1}\right)\\
&=&\left(ab_1a^{-1}b_1^{-1}\right)\left(b_1ab_2a^{-1}b_2^{-1}b_1^{-1}\right)\\
&=&ab_1b_2a^{-1}b_2^{-1}b_1^{-1}\\
&=&[a,b_1b_2],\\
&&[a,b_1b_2\cdots{}b_{n-1}]\left(b_1b_2\cdots{}b_{n-1}[a,b_n]b_{n-1}^{-1}\cdots{}b_2^{-1}b_1^{-1}\right)\\
&=&\left(ab_1b_2\cdots{}b_{n-1}a^{-1}b_{n-1}^{-1}\cdots{}b_2^{-1}b_1^{-1}\right)\left(b_1b_2\cdots{}b_{n-1}ab_na^{-1}b_n^{-1}b_{n-1}^{-1}\cdots{}b_2^{-1}b_1^{-1}\right)\\
&=&ab_1b_2\cdots{}b_na^{-1}b_n^{-1}\cdots{}b_2^{-1}b_1^{-1}\\
&=&[a,b_1b_2\cdots{}b_n].
\end{eqnarray*}
By induction on $n$, we obtain the claim.
\end{proof}

\begin{proof}[Proof of Lemma~\ref{gen-T^2I-1}]
By Theorem~\ref{main-1}, we have the short exact sequence
$$1\to\T^2(N_g)\I(N_g)\to\M_2(N_g)\to\left(\Z/2\Z\right)^{\tbinom{g-1}{2}+\epsilon}\to1,$$
where $\epsilon=0$ if $g$ is odd and  $\epsilon=1$ if $g$ is even.
In order to give a generating set for $\T^2(N_g)\I(N_g)$, we use the Reidemeister-Schreier method (for details, for instance see \cite{Jo}).

By Lemma~\ref{G/N-1}, it follows that $\left(\Z/2\Z\right)^{\tbinom{g-1}{2}+\epsilon}$ is equal to $\YLO$ as a set.
We note that $\YLO$ is a Schreier transversal for $\T^2(N_g)\I(N_g)$ in $\M_2(N_g)$, that is, for any $Y_{\alpha_{i_1},\alpha_{i_1,j_1}}Y_{\alpha_{i_2},\alpha_{i_2,j_2}}\cdots{}Y_{\alpha_{i_m},\alpha_{i_m,j_m}}\in\YLO$ with $m\geq1$, $Y_{\alpha_{i_1},\alpha_{i_1,j_1}}Y_{\alpha_{i_2},\alpha_{i_2,j_2}}\cdots{}Y_{\alpha_{i_{m-1}},\alpha_{i_{m-1},j_{m-1}}}$ is in $\YLO$.
For any $x\in\M_2(N_g)$, let $\overline{x}$ be the element of $\YLO=\left(\Z/2\Z\right)^{\tbinom{g-1}{2}+\epsilon}$ which is sent from $x$ by the natural projection.
Then $\T^2(N_g)\I(N_g)$ is generated by
$$G=\left\{\left.fx^{\pm1}\overline{fx^{\pm1}}^{-1}\right|x\in{X},f\in\YLO,fx^{\pm1}\neq\overline{fx^{\pm1}}\right\},$$
where $X$ is the generating set for $\M_2(N_g)$ in Theorem~\ref{level2-gen}.
For any $fx^{\pm1}\overline{fx^{\pm1}}^{-1}\in{G}$, we check that $fx^{\pm1}\overline{fx^{\pm1}}^{-1}$ is a product of elements in (1)-(4) of Lemma~\ref{gen-T^2I-1} and their inverses.
Note that $Y_{\alpha_j,\alpha_{i,j}}^{-1}Y_{\alpha_i,\alpha_{i,j}}=Y_{\alpha_j,\alpha_{i,j}}Y_{\alpha_i,\alpha_{i,j}}^{-1}=t_{\alpha_{i,j}}^2$ and $Y_{\alpha_i,\alpha_{i,j}}^2=Y_{\alpha_j,\alpha_{i,j}}^2=t_{\beta_{i,j}}$.

Let $x=Y_{\alpha_i,\alpha_{i,j}}$.
If $f$ includes $x$, there are $f_1$, $f_2\in\YLO$ such that $f=f_1xf_2$.
Then we see
\begin{eqnarray*}
fx\overline{fx}^{-1}
&=&f_1xf_2x\overline{f_1xf_2x}^{-1}\\
&=&f_1xf_2x(f_1f_2)^{-1}\\
&=&f_1[x,f_2]f_1^{-1}\cdot{}f_1f_2x^2f_2^{-1}f_1^{-1},\\
fx^{-1}\overline{fx^{-1}}^{-1}
&=&f_1xf_2x^{-1}\overline{f_1xf_2x^{-1}}^{-1}\\
&=&f_1xf_2x^{-1}(f_1f_2)^{-1}\\
&=&f_1[x,f_2]f_1^{-1}.
\end{eqnarray*}
If $f$ does not include $x$, there are $f_1$, $f_2\in\YLO$ such that $f=f_1f_2$ and $f_1xf_2\in\YLO$.
Then we see
\begin{eqnarray*}
fx\overline{fx}^{-1}
&=&f_1f_2x\overline{f_1f_2x}^{-1}\\
&=&f_1f_2x(f_1xf_2)^{-1}\\
&=&f_1[x,f_2]^{-1}f_1^{-1},\\
fx^{-1}\overline{fx^{-1}}^{-1}
&=&f_1f_2x^{-1}\overline{f_1f_2x^{-1}}^{-1}\\
&=&f_1f_2x^{-1}(f_1xf_2)^{-1}\\
&=&fx^{-2}f^{-1}\cdot{}f_1[x,f_2]^{-1}f_1^{-1}.
\end{eqnarray*}
By Lemma~\ref{[a,b]}, we have that $fx^{\pm1}\overline{fx^{\pm1}}^{-1}$ is a product of elements in (1)-(4) of Lemma~\ref{gen-T^2I-1} and their inverses.
In addition, let $x=Y_{\alpha_j,\alpha_{i,j}}$ and $x^\prime=Y_{\alpha_i,\alpha_{i,j}}$.
Then we see
\begin{eqnarray*}
fx\overline{fx}^{-1}&=&fx\overline{fx^\prime}^{-1}=fx(x^\prime)^{-1}f^{-1}\cdot{}fx^\prime\overline{fx^\prime}^{-1},\\
fx^{-1}\overline{fx^{-1}}^{-1}&=&fx^{-1}\overline{fx}^{-1}=fx^{-2}f^{-1}\cdot{}fx\overline{fx}^{-1}.
\end{eqnarray*}
Moreover, let $x=t_{\alpha_{1,j,k,l}}^2$.
Then we see
$$fx^{\pm1}\overline{fx^{\pm1}}^{-1}=fx^{\pm1}f^{-1}.$$
Therefore, for any $fx^{\pm1}\overline{fx^{\pm1}}^{-1}\in{G}$, we have that $fx^{\pm1}\overline{fx^{\pm1}}^{-1}$ is a product of elements in (1)-(4) of Lemma~\ref{gen-T^2I-1} and their inverses.

Thus we finish the proof.
\end{proof}

\begin{rem}\label{reduced-gen}
Put $f=Y_{\alpha_{i_1},\alpha_{i_1,j_1}}Y_{\alpha_{i_2},\alpha_{i_2,j_2}}\cdots{}Y_{\alpha_{i_m},\alpha_{i_m,j_m}}\in\YLO$.
We notice that the generator (4) in Lemma~\ref{gen-T^2I-1} can be reduced a form satisfying $(i_m,j_m)<(i,j)<(k,l)$, since in the above proof $f_1[x,f_2]^{-1}f_1^{-1}$ is a product of these reduced generators.
\end{rem}

\begin{proof}[Proof of Theorem~\ref{gen-T^2I-2}]
It is clear that the generators (1), (2) and (3) in Lemma~\ref{gen-T^2I-1} are equal to the generators (1), (2) and (3) in Theorem~\ref{gen-T^2I-2} respectively.
We show that the reduced generator (4) in Remark~\ref{reduced-gen} is obtained from the generators (1) and (2) in Theorem~\ref{gen-T^2I-2}.

For any reduced generator $f[Y_{\alpha_{i},\alpha_{i,j}},Y_{\alpha_{k},\alpha_{k,l}}]f^{-1}$, we see
\begin{eqnarray*}
f[Y_{\alpha_{i},\alpha_{i,j}},Y_{\alpha_{k},\alpha_{k,l}}]f^{-1}
&=&
f(Y_{\alpha_{i},\alpha_{i,j}}Y_{\alpha_{k},\alpha_{k,l}})^2Y_{\alpha_{k},\alpha_{k,l}}^{-2}(Y_{\alpha_{k},\alpha_{k,l}}Y_{\alpha_{i},\alpha_{i,j}}^{-2}Y_{\alpha_{k},\alpha_{k,l}}^{-1})f^{-1}\\
&=&
f(Y_{\alpha_{i},\alpha_{i,j}}Y_{\alpha_{k},\alpha_{k,l}})^2f^{-1}\cdot{}ft_{\beta_{k,l}}^{-1}f^{-1}\\
&&
f(Y_{\alpha_{k},\alpha_{k,l}}t_{\beta_{i,j}}^{-1}Y_{\alpha_{k},\alpha_{k,l}}^{-1})f^{-1}.
\end{eqnarray*}
Since $f$ and $fY_{\alpha_{k},\alpha_{k,l}}$ are in $\YLO$, $ft_{\beta_{k,l}}f^{-1}$ and $f(Y_{\alpha_{k},\alpha_{k,l}}t_{\beta_{i,j}}Y_{\alpha_{k},\alpha_{k,l}}^{-1})f^{-1}$ are both generators (2) in Theorem~\ref{gen-T^2I-2}.
Hence we check that either $f[Y_{\alpha_{i},\alpha_{i,j}},Y_{\alpha_{k},\alpha_{k,l}}]f^{-1}$ or $f(Y_{\alpha_{i},\alpha_{i,j}}Y_{\alpha_{k},\alpha_{k,l}})^2f^{-1}$ is obtained from the generators (1) and (2) in Theorem~\ref{gen-T^2I-2}.

Let $i=k$.
Then we calculate
\begin{eqnarray*}
(Y_{\alpha_{i},\alpha_{i,j}}Y_{\alpha_{k},\alpha_{k,l}})^2
&=&
\left(t_{\alpha_{j,l}}t_{Y_{\alpha_i,\alpha_{i,j}}^{-1}(\alpha_{j,l})}^{-1}\right)^2\\
&=&
t_{\alpha_{j,l}}^2(Y_{\alpha_i,\alpha_{i,j}}^{-1}t_{\alpha_{j,l}}^{-2}Y_{\alpha_i,\alpha_{i,j}})\\
&=&
t_{\alpha_{j,l}}^2Y_{\alpha_i,\alpha_{i,j}}^{-2}(Y_{\alpha_i,\alpha_{i,j}}t_{\alpha_{j,l}}^{-2}Y_{\alpha_i,\alpha_{i,j}}^{-1})Y_{\alpha_i,\alpha_{i,j}}^2\\
&=&
t_{\alpha_{j,l}}^2t_{\beta_{i,j}}^{-1}(Y_{\alpha_i,\alpha_{i,j}}t_{\alpha_{j,l}}^{-2}Y_{\alpha_i,\alpha_{i,j}}^{-1})t_{\beta_{i,j}}
\end{eqnarray*}
(see Figure~\ref{crosscap-slide-rel-3}~(a)).
Since $f$ and $fY_{\alpha_i,\alpha_{i,j}}$ are in $\YLO$, $f(Y_{\alpha_{i},\alpha_{i,j}}Y_{\alpha_{k},\alpha_{k,l}})^2f^{-1}$ is obtained from the generators (1) and (2) in Theorem~\ref{gen-T^2I-2}.
Let $j=l$.
Then we calculate
\begin{eqnarray*}
[Y_{\alpha_{i},\alpha_{i,j}},Y_{\alpha_{k},\alpha_{k,l}}]
&=&
Y_{\alpha_{i},\alpha_{i,j}}Y_{\alpha_{k},\alpha_{k,j}}Y_{\alpha_{i},\alpha_{i,j}}^{-1}Y_{\alpha_{k},\alpha_{k,j}}^{-1}\\
&=&
Y_{Y_{\alpha_{i},\alpha_{i,j}}(\alpha_{k}),Y_{\alpha_{i},\alpha_{i,j}}(\alpha_{k,j})}Y_{\alpha_{k},\alpha_{k,j}}^{-1}\\
&=&
Y_{\alpha_{k},\alpha_{k,j}}Y_{\alpha_{k},\alpha_{i,k}}^2Y_{\alpha_{k},\alpha_{k,j}}^{-1}\\
&=&
Y_{\alpha_{k},\alpha_{k,j}}t_{\beta_{i,k}}Y_{\alpha_{k},\alpha_{k,j}}^{-1}
\end{eqnarray*}
(see Figure~\ref{crosscap-slide-rel-3}~(b)).
Since $fY_{\alpha_k,\alpha_{k,j}}$ is in $\YLO$, $f[Y_{\alpha_{i},\alpha_{i,j}},Y_{\alpha_{k},\alpha_{k,l}}]f^{-1}$ is obtained from the generators (1) and (2) in Theorem~\ref{gen-T^2I-2}.
Let $j=k$.
Then we calculate
\begin{eqnarray*}
(Y_{\alpha_{i},\alpha_{i,j}}Y_{\alpha_{k},\alpha_{k,l}})^2
&=&
Y_{\alpha_{i},\alpha_{i,j}}Y_{\alpha_{j}\alpha_{j,l}}Y_{\alpha_{i},\alpha_{i,j}}Y_{\alpha_{j},\alpha_{i,j}}^{-1}Y_{\alpha_{j}\alpha_{j,l}}^{-1}Y_{\alpha_{i},\alpha_{i,j}}^{-1}\\
&&
Y_{\alpha_{i},\alpha_{i,j}}Y_{\alpha_{j},\alpha_{i,j}}^{-1}(Y_{\alpha_{j},\alpha_{i,j}}Y_{\alpha_{j},\alpha_{j,l}})^2\\
&=&
Y_{\alpha_{i},\alpha_{i,j}}Y_{\alpha_{j}\alpha_{j,l}}t_{\alpha_{i,j}}^2Y_{\alpha_{j}\alpha_{j,l}}^{-1}Y_{\alpha_{i},\alpha_{i,j}}^{-1}\\
&&
t_{\alpha_{i,j}}^2\left(t_{Y_{\alpha_j,\alpha_{j,l}}^{-2}(\alpha_{i,l})}t_{Y_{\alpha_j,\alpha_{j,l}}^{-1}(\alpha_{i,l})}^{-1}\right)^2\\
&=&
Y_{\alpha_{i},\alpha_{i,j}}Y_{\alpha_{j}\alpha_{j,l}}t_{\alpha_{i,j}}^2Y_{\alpha_{j}\alpha_{j,l}}^{-1}Y_{\alpha_{i},\alpha_{i,j}}^{-1}\\
&&
t_{\alpha_{i,j}}^2(Y_{\alpha_j,\alpha_{j,l}}^{-2}t_{\alpha_{i,l}}^2Y_{\alpha_j,\alpha_{j,l}}^2\cdot{}Y_{\alpha_j,\alpha_{j,l}}^{-1}t_{\alpha_{i,l}}^{-2}Y_{\alpha_j,\alpha_{j,l}})\\
&=&
Y_{\alpha_{i},\alpha_{i,j}}Y_{\alpha_{j}\alpha_{j,l}}t_{\alpha_{i,j}}^2Y_{\alpha_{j}\alpha_{j,l}}^{-1}Y_{\alpha_{i},\alpha_{i,j}}^{-1}\\
&&
t_{\alpha_{i,j}}^2t_{\beta_{j,l}}^{-1}t_{\alpha_{i,l}}^2t_{\beta_{j,l}}\cdot{}Y_{\alpha_j,\alpha_{j,l}}t_{\beta_{j,l}}^{-1}t_{\alpha_{i,l}}^{-2}t_{\beta_{j,l}}Y_{\alpha_j,\alpha_{j,l}}^{-1}
\end{eqnarray*}
(see Figure~\ref{crosscap-slide-rel-3}~(c)).
Since $f$, $fY_{\alpha_j,\alpha_{j,l}}$ and $fY_{\alpha_i,\alpha_{i,j}}Y_{\alpha_j,\alpha_{j,l}}$ are in $\YLO$, $f(Y_{\alpha_{i},\alpha_{i,j}}Y_{\alpha_{k},\alpha_{k,l}})^2f^{-1}$ is obtained from the generators (1) and (2) in Theorem~\ref{gen-T^2I-2}.
Let $i<k<j<l$.
Then we calculate
\begin{eqnarray*}
[Y_{\alpha_{i},\alpha_{i,j}},Y_{\alpha_{k},\alpha_{k,l}}]
&=&
Y_{\alpha_{i},\alpha_{i,j}}Y_{\alpha_{k},\alpha_{k,l}}Y_{\alpha_{i},\alpha_{i,j}}^{-1}Y_{\alpha_{k},\alpha_{k,l}}^{-1}\\
&=&
Y_{\alpha_{i},\alpha_{i,j}}Y_{Y_{\alpha_{k},\alpha_{k,l}}(\alpha_{i}),Y_{\alpha_{k},\alpha_{k,l}}(\alpha_{i,j})}^{-1}\\
&=&
Y_{\alpha_{i},\alpha_{i,j}}Y_{\alpha_{i},\alpha_{i,l}}^{-1}Y_{\alpha_{i},\alpha_{i,k}}^2Y_{\alpha_{i},\alpha_{i,l}}Y_{\alpha_{i},\alpha_{i,k}}^2Y_{\alpha_{i},\alpha_{i,j}}^{-1}\\
&&
Y_{\alpha_{i},\alpha_{i,k}}^{-2}Y_{\alpha_{i},\alpha_{i,l}}^{-1}Y_{\alpha_{i},\alpha_{i,k}}^{-2}Y_{\alpha_{i},\alpha_{i,l}}\\
&=&
Y_{\alpha_{i},\alpha_{i,j}}(t_{\beta_{i,l}}^{-1}Y_{\alpha_{i},\alpha_{i,l}}t_{\beta_{i,k}}Y_{\alpha_{i},\alpha_{i,l}}^{-1}t_{\beta_{i,l}}t_{\beta_{i,k}})Y_{\alpha_{i},\alpha_{i,j}}^{-1}\\
&&
t_{\beta_{i,k}}^{-1}t_{\beta_{i,l}}^{-1}\cdot{}Y_{\alpha_{i},\alpha_{i,l}}t_{\beta_{i,k}}^{-1}Y_{\alpha_{i},\alpha_{i,l}}^{-1}\cdot{}t_{\beta_{i,l}}
\end{eqnarray*}
(see Figure~\ref{crosscap-slide-rel-3}~(d)).
Since $f$, $fY_{\alpha_i,\alpha_{i,j}}$, $fY_{\alpha_i,\alpha_{i,l}}$ and $fY_{\alpha_i,\alpha_{i,j}}Y_{\alpha_i,\alpha_{i,l}}$ are in $\YLO$, $f[Y_{\alpha_{i},\alpha_{i,j}},Y_{\alpha_{k},\alpha_{k,l}}]f^{-1}$ is obtained from the generators (1) and (2) in Theorem~\ref{gen-T^2I-2}.
For the other cases, we have that $f[Y_{\alpha_{i},\alpha_{i,j}},Y_{\alpha_{k},\alpha_{k,l}}]f^{-1}=1$.

\begin{figure}[htbp]
\subfigure[$Y_{\alpha_{i},\alpha_{i,j}}Y_{\alpha_{i},\alpha_{i,l}}=t_{\alpha_{j,l}}t_{Y_{\alpha_i,\alpha_{i,j}}^{-1}(\alpha_{j,l})}^{-1}$.]{\includegraphics{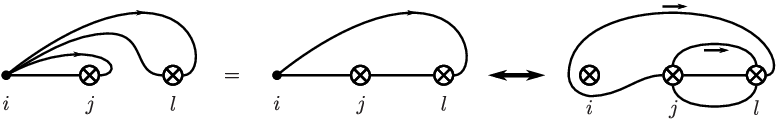}}
\subfigure[$Y_{Y_{\alpha_{i},\alpha_{i,j}}(\alpha_{k}),Y_{\alpha_{i},\alpha_{i,j}}(\alpha_{k,j})}=Y_{\alpha_{k},\alpha_{k,j}}Y_{\alpha_{k},\alpha_{i,k}}^2$.]{\includegraphics{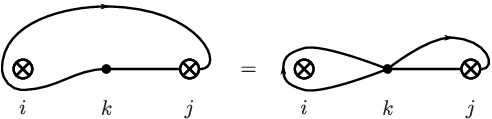}}
\subfigure[$Y_{\alpha_{j},\alpha_{i,j}}Y_{\alpha_{j},\alpha_{j,l}}=t_{Y_{\alpha_j,\alpha_{j,l}}^{-2}(\alpha_{i,l})}t_{Y_{\alpha_j,\alpha_{j,l}}^{-1}(\alpha_{i,l})}^{-1}$.]{\includegraphics{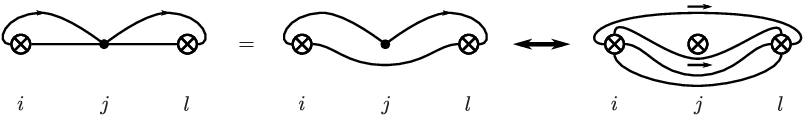}}
\subfigure[$Y_{Y_{\alpha_{k},\alpha_{k,l}}(\alpha_{i}),Y_{\alpha_{k},\alpha_{k,l}}(\alpha_{i,j})}^{-1}=Y_{\alpha_{i},\alpha_{i,l}}^{-1}Y_{\alpha_{i},\alpha_{i,k}}^2Y_{\alpha_{i},\alpha_{i,l}}Y_{\alpha_{i},\alpha_{i,k}}^2Y_{\alpha_{i},\alpha_{i,j}}^{-1}Y_{\alpha_{i},\alpha_{i,k}}^{-2}Y_{\alpha_{i},\alpha_{i,l}}^{-1}Y_{\alpha_{i},\alpha_{i,k}}^{-2}Y_{\alpha_{i},\alpha_{i,l}}$.]{\includegraphics{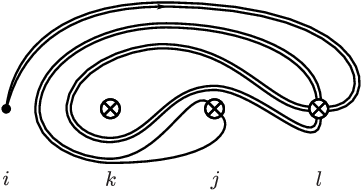}}
\caption{}\label{crosscap-slide-rel-3}
\end{figure}

Thus we complete the proof.
\end{proof}

\appendix
\section{The case where $g$ is odd}\label{T^2I}

In this appendix, we prove the following.

\begin{thm}\label{main-2}
$\I(N_g)\subset\T^2(N_g)$ if $g\geq3$ is odd.
\end{thm}

\begin{proof}
Recall Subsection~\ref{M_2-I}, especially Theorem~\ref{torelli-gen}.
Since $\I(N_3)$ is trivial, it is clear that $\I(N_3)\subset\T^2(N_3)$.
Hence suppose that $g\geq5$.
By Lemma~\ref{chain-T^2}, it follows that $t_{\beta_1}t_{\beta_2}^{-1}(=t_{\beta_1}t_{\beta_2}\cdot{}t_{\beta_2}^{-2})$ is in $\T^2(N_g)$.
Note that $N_g\setminus\alpha$ is a joint union of $N_2$ with a boundary and $N_{g-2}$ with a boundary.
Let $\delta$ be a simple closed curve such that $N_g\setminus\left(\alpha\cup\delta\right)$ is a joint union of a M\"obius band, $N_2$ with a boundary and $\Sigma_{\frac{g-3}{2}}$ with two boundaries, as shown in Figure~\ref{alpha-delta}.
By Lemma~\ref{chain-T^2}, it follows that $t_\alpha{}t_\delta$ is in $\T^2(N_g)$.
Since $t_\delta=1$, we have that $t_\alpha$ is in $\T^2(N_g)$.
Therefore we conclude that $\I(N_g)\subset\T^2(N_g)$.
Thus we finish the proof.

\begin{figure}[htbp]
\includegraphics{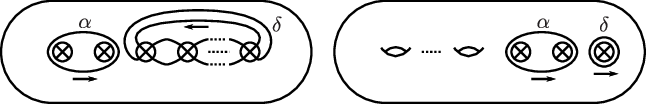}
\caption{Simple closed curves $\alpha$ and $\delta$.
The figure on the left side is a connected sum of $g$ real projective planes.
The figure on the right side is a connected sum of $\Sigma_{\frac{g-3}{2}}$ and $3$ real projective planes.}\label{alpha-delta}
\end{figure}

\end{proof}

By Theorems~\ref{main-1} and $\ref{main-2}$, we have the following.

\begin{cor}
For odd $g\geq3$, the quotient of $\M_2(N_g)$ by $\T^2(N_g)$ is isomorphic to $\left(\Z/2\Z\right)^{\tbinom{g-1}{2}}$, and the quotient of $\T_2(N_g)$ by $\T^2(N_g)$ is isomorphic to $\left(\Z/2\Z\right)^{\tbinom{g-1}{2}-1}$.
\end{cor}

\begin{rem}
Any Dehn twist about a separating simple closed curve is in $\I(N_g)$ for any genus $g$, and so in $\T^2(N_g)$ if $g$ is odd.
\end{rem}

\begin{rem}
We do not know whether or not $\I(N_g)\subset\T^2(N_g)$ if $g$ is even.
\end{rem}


\end{document}